\documentclass[10pt]{amsart}

\usepackage{color}
\usepackage[colorlinks=true,linkcolor=blue]{hyperref} 

\usepackage{graphics}
\usepackage{epsfig} 
\usepackage[colorlinks=true,linkcolor=blue]{hyperref} 
\newtheorem{theor}{{\bf Theorem}}[section]
\newtheorem{lemma}[theor]{{\bf Lemma}}           
\newtheorem{coro}[theor]{{\bf Corollary}}

\newtheorem{definition}[theor]{{\bf Definition}}
\newtheorem{remark}{{\bf Remark}}[section]
\numberwithin{equation}{section}

\newcommand{\refq}[1]{~(\ref{#1})}

\newcommand{\egal}{\Longleftrightarrow}
\newcommand{\sauf}{\setminus}

\newcommand{\epsl}{\varepsilon}

\newcommand{\bfepsl}{\boldsymbol{\epsl}} 
\newcommand{\bfmu}{\boldsymbol{\mu}}

\newcommand{\calA}{\mathcal A}

\newcommand{\calD}{\mathcal D}
\newcommand{\calE}{\mathcal E}
\newcommand{\calF}{\mathcal F}

\newcommand{\calH}{\mathcal H}

\newcommand{\calM}{\mathcal M}

\newcommand{\calV}{\mathcal V}

\newcommand{\HDp}{H^{D_{\rm a}}}

\newcommand{\rotdisc}{{\nabla_{\rm disc}\wedge}} 

\newcommand{\partOm}{\partial\Om}
\newcommand{\partanOm}{\partial_{\rm tan}}
\newcommand{\tanOm}{{_{{\rm tan}}}}

\newcommand{\Hadm}{\calH_{\partOm,{\rm adm}}}
\newcommand{\HadmE}{\calH_{\Gamma,{\rm adm}}}
\newcommand{\HadmEp}{\calH_{\partOm\sauf\partOm_1^+,{\rm adm}}}
\newcommand{\HadmEm}{\calH_{\partOm_1^-,{\rm adm}}}

\newcommand{\Om}{\Omega}

\newcommand{\om}{\omega}

\newcommand{\CC}{\mathbf C}
\newcommand{\R}{\mathbf R}
\newcommand{\Z}{\mathbf Z}

\title{Inverse problem for the discrete Maxwell Equations in a bounded paving}
\author{Olivier Poisson$^{\small 1}$}

\begin{document}
\baselineskip 14pt

\maketitle
%
{\footnotesize
 \centerline{$^{\small 1}$ Aix-Marseille Universit\'{e}, France}
} 

\begin{abstract}
 We consider the discrete anisotropic Maxwell operator $D_{\rm a} H_0$ on
 a bounded paving $\Om\subset\Z^3$, where $H_0$ denotes discrete isotropic
 Maxwell operator and $D_{\rm a}$ a diagonal operator of multiplication
 containing information about the anisotropy of the medium inside $\Om$.
 Letting a complex number $\lambda\neq 0$ such the Dirichlet-to-Neumann
 operator $\Lambda(D_{\rm a})$ associated with the system $D_{\rm a} H_0 u =\lambda u$
 on $\Om$ admits a unique solution, we show that knowing $\Lambda(D_{\rm a})$
 is sufficient to determine $D_{\rm a}$ by a reconstruction procedure for $D_{\rm a}$.
\end{abstract}



\section{Introduction}
\subsection{Presentation}
\subsubsection*{The formal discrete Maxwell operator}
 Let $\bfepsl$ and $\bfmu$ be the permittivity and the permeability in the ambient space $\Z^3$, 
 which are the $3\times 3$ constant diagonal matrices with diagonal elements,
 $\epsl_1, \epsl_2, \epsl_3\in (0,\infty)$ for $\bfepsl$, and $\mu_1,\mu_2,\mu_3 \in (0,\infty)$
 for $\bfmu$.
 Letting $\Z^d=\{n=(n_1,\ldots,n_d);\; n_j\in\Z\}$, $d\ge 1$, the square lattice,
 setting $(e_1,\ldots,e_d)$ the usual euclidian basis of $\R^d$, and denoting by $\calF(E;F)$
 (respect., by $\calA(E;F)$) the set of functions (respect., applications) from a set $E$
 into another set $F$, we put, for $u=(u^E,u^H)\in \calF(\Z^3;\CC^6)$ with $u^E,u^H \in \calF(\Z^3;\CC^3)$,
 the isotropic discrete Maxwell operator:
\begin{equation}
\label{def.H0u}
   H_0 u = (M u^H, -M u^E),
\end{equation}
 where $M$ denotes the discrete rotational operator defined by
\begin{equation}
\label{def.hatM}
   M v = (2i)^{-1}  \rotdisc v  :=
  (2i)^{-1} (- D_3 v_2 + D_2 v_3, D_3 v_1 - D_1 v_3, - D_2 v_1 + D_1 v_2),
\end{equation}
 with $v=(v_1,v_2,v_3) \in \calF(Z^3;\CC^3)$,
 and $D_j$ denotes the discrete difference operator defined by
$$ D_j u\; (n) :=   u(n+e_j)- u(n-e_j) , \quad j\in [[1,d]].$$ 
%

%
 %
 
  Let $D_{\rm a}:\: \Z^3 \to \R^6$ be a function such that $D_{\rm a}(n)$ is a real positive diagonal
 of the form
\[
   D_{\rm a} = \left(\begin{array}{cc} \bfepsl_p & 0_{3\times3} \\ 0_{3\times3} &\bfmu_p
  \end{array}\right).\qedhere
\]
 We assume that $(\bfepsl_p,\bfmu_p)=(\bfepsl_0,\bfmu_0)$ is constant outside a compact set,
 $\Om$. 
 We put $D_0= D_{\rm a}$ when $(\bfepsl_p,\bfmu_p)=(\bfepsl_0,\bfmu_0)$.
 The anisotropic Maxwell operator is defined by
\[  \HDp =   D_{\rm a}  H_0.\]

\subsubsection*{Geometry on the lattice}
 We define the following relation of equivalence on $\Z^d$, $d\ge 1$:
$$ n\sim m \egal |m-n|=1, \quad n,m\in \Z^d , $$ 
 so we have $\{m;\; n\sim m\} = \cup_{j=1}^d \{n\pm e_j\}$.
 For $n\in\Z^d$ and $\Om\subset\Z^d$, we write $n\sim \Om$ means that there exists $m\in \Om$
 such that $n\sim m$. 
 We also define 
$$ \Om' = \{n\in \Z^d ; \; n\sim \Om \} , \quad \partOm = \Om'\sauf \Om. $$
 We call $\partOm$ the boundary of the open set $\Om$.

 From now we consider the case where $\Om$ is a paving, i.e. a bounded set of the form
 $\Om=\Pi_{j=1}^d[[a_j,b_j]]\subset\Z^d$ with $a_j, b_j\in\Z$.
 Letting $n\in\partOm$, there exists then exactly one point $m=m_\Om(n)\in\Om$ such that
 $n\sim m$. Hence, we define the (outgoing) normal vector field at the boundary $\partOm$:
$$ \nu(n)=n-m_\Om(n) \in \{\pm e_j;\; j\in [[1,d]]\},\quad n\in\partOm.$$ 
 From now we consider the case $d=3$.
 Setting
\begin{eqnarray*}
 \partOm_k^+ &:=& \{n\in \Z^3; \; \quad n_k = b_k+1, \; n_j\in [[a_j,b_j]] 
  \quad  \forall j\neq k\} ,\\
 \partOm_k^- &:=& \{n\in \Z^3; \; \quad n_k = a_k-1, \; n_j\in [[a_j,b_j]]
  \quad  \forall j\neq k\} ,\\
 \partOm_k &:=&  \partOm_k^+ \cup  \partOm_k^- =
  \{n\in \partOm; \; n_j\in [[a_j,b_j]] , \quad  \forall j\neq k\} ,
\end{eqnarray*}
 we then can write the boundary as the partitions:
$$ \partOm = \cup_j \partOm_j= \cup_\pm \cup_j \partOm_j^\pm ,$$
 and we have, for $n\in\partOm_k$, $\nu(n)=e_k$ if $n_k=b_k$ and $\nu(n)=-e_k$ if $n_k=a_k$.
 We put, for $u\in \calA(\Om';\CC^{3})$,
$$ u\tanOm(n) := u(n) - (u(n)\cdot\nu(n))\nu(n) , \quad n\in \partOm ,$$ 
$$ \partanOm u(n) = (u(n)-u(m_\Om(n))) \wedge \nu(n), \quad n\in \partOm,$$ 
 and, for $l\ge 1$, $u=(u^{(1)},\ldots,u^{(l)})\in \calA(\Om';\CC^{3l})$,
$$
  u\tanOm=(u^{(1)}\tanOm,\ldots,u^{(l)}\tanOm), \quad \partanOm u(n)
   = (\partanOm u^{(1)},\ldots,\partanOm u^{(l)}).
$$
 In addition, writing $u^{(s)}=(u_1^{(s)},u_2^{(s)},u_3^{(s)})\in \calA(\Om';\CC^3)$,
 we see that $u^{(s)}\tanOm|_{\partOm_j}$ does not depend on $u^{(s)}_j$.
%
 Hence, setting 
\begin{center}
 $\calH(\Om;\CC^{3l}):= \{u=(u^{(1)},\ldots,u^{(l)})$; $u^{(s)}\in \calA(\Om\cup\partOm_k;\CC^l)$, $k\neq j$,
 $s\in [[1,l]]\}$,
\end{center}
 and the admissible set of data on $\partOm$ as
\begin{equation}
\label{def.Hadm}
 \Hadm:=\{f=(f_1,f_2,f_3)\in \calA(\partOm;\CC^3);\;  f(n)\cdot \nu(n)=0 ,\: \forall n\in\partOm\},
\end{equation}
 we have $u\tanOm,\partanOm u\in (\Hadm)^{l}$ for all $u=\calH(\Om;\CC^{3l})$.
 
\subsubsection*{The discrete Maxwell operator on the paving $\Om$} 
 For simplicity, we assume that $\Om=[[1,R_1]]\times[[1,R_2]]\times[[1,R_3]]$
 from now.
 Letting $\lambda\in\CC$, we consider on $\calH(\Om;\CC^6)$ the homogeneous
 Maxwell system:
\begin{equation}
\label{eq.HDu=lambdau}
  \HDp u -\lambda u = 0 ,
\end{equation}
 with the boundary condition
\begin{equation}
\label{eq.u=f}
 u_\tanOm = f \quad {\rm in} \quad \partOm,
\end{equation}
 where $f=(f^E,f^H) \in (\Hadm)^2$ is given.
 Obviously, the system\refq{eq.HDu=lambdau}-\refq{eq.u=f} is finite linear square
 with $6R_1R_2R_3+4(R_1R_2+R_2R_3+R_3R_1)$ unknowns and scalar equations.
 (Observe that $2(R_1R_2+R_2R_3+R_3R_1)=\dim\Hadm$.)
 We denote by $\calE(D_{\rm a})$ the finite set of eigenvalues of the above system.
 For $\lambda\not\in \calE(D_{\rm a})$, System\refq{eq.HDu=lambdau} admits a unique solution,

\subsubsection*{The Calderon inverse problem}
 For $\lambda\not\in \calE(D_{\rm a})$, we then define the following operator on $\calH(\Om;\CC^6)$:
$$ \Lambda(D_{\rm a}) \: : \quad  f \mapsto \partanOm u , $$
 where $u$ is the solution of\refq{eq.HDu=lambdau}.
 This is the Dirichlet-to-Neumann ("D-N") operator.

\medskip
 The Calderon inverse problem for the discrete Maxwell anisotropic operator
 can be formulated as follows:
\begin{center}
Q: {\em knowing the background $(\bfepsl_0,\bfmu_0)$, does the knowledge of $\Lambda(D_{\rm a})$
 determines $D_{\rm a}$?}
\end{center} 

\subsection{Main result}
 In this article we prove that the response to question Q is yes.
 Let us give more precisions. 
\begin{theor}
\label{th.1}
 Let two anisotropic operators $H^{D_{\rm a}}= D_{\rm a} H_0$ and $H^{\tilde D_{\rm a}}= \tilde D_{\rm a} H_0$
 with common background $(\bfepsl_0,\bfmu_0)$. 
 Assume that $\lambda \not\in \calE(D_{\rm a})\cup \calE(\tilde D_{\rm a})$.
 If $\Lambda(D_{\rm a})=\Lambda(\tilde D_{\rm a})$ we then have $D_{\rm a}=\tilde D_{\rm a}$.

 In addition, starting from the knowledge of $\Lambda(D_{\rm a})$, we give an explicit
 reconstruction of $D_{\rm a}$.
\end{theor}
\begin{remark}
\begin{itemize}
\item If $\lambda\in \calE(D_{\rm a})$ the D-N operator is undefined, but it is possible
 to reformulate the Calderon problem in terms of the set of Cauchy data
 $(u_\tanOm,\partanOm u)$ where $u$ satisfies \refq{eq.HDu=lambdau}.
 If $\lambda=0$, the Cauchy data do not depend on $D_{\rm a}$
 since\refq{eq.HDu=lambdau} is equivalent to the equation $H_0u=0$.
 However, in applications, $\lambda$ is usually non-zero.\\
 From now we assume that $\lambda\neq \calE(D_{\rm a})$.
\item The proof shows that the above result extends to the case where, for all $n$,
 $D_{\rm a}(n)$ is a complex-valued invertible diagonal matrix.
\item Our proof is constructive:
 we compute $D_{\rm a}$ from the construction of solutions of\refq{eq.HDu=lambdau}
 which depend only on $\Lambda(D_{\rm a})$.
\item Actually, our proof treats the case of partial data since, in the below
 reconstruction of $D_{\rm a}$, we impose $u_\tanOm=0$ on $\partOm_3\cup\partOm_2^-$.
 It means that $D_{\rm a}$ can be reconstructed from the following D-N operator:
$$
  (\Hadm)^2 \ni f  \mapsto \partanOm u , 
$$
 where $u\in \calH(\Om;\CC^6)$ is the unique solution of\refq{eq.HDu=lambdau}
 such that $u_{\tanOm}=f$ on $\partOm_1\cup  \partOm_2^+$
 and $u_{\tanOm}=0$ on $\partOm_2^-\cup  \partOm_3$.\\
 In addition, it is possible to attribute for $u_\tanOm$ on $\partOm_2^-\cup \partOm_3$ 
 any other fixed values.
\end{itemize} 

\end{remark} 

\subsection{Motivation and Related works}
 In the continuous case, the Calderon inverse problem related to the Maxwell system
 in a isotropic medium with  assuming permeability $\bfmu=1$ was solved in \cite{Co-Pa}. 
 In \cite{Ch-Ol} and \cite{Ol-So} the proof of uniqueness was extended to
 unknown permeability, permissibility and conductivity by using generalized
 Sommerfeld potentials.
 In \cite{Ke-Sa}, the proof of uniqueness was extended to the large framework
 of admissible Riemannian manifolds, so, it works for anisotropic Maxwell systems
 in the Euclidean space.\\
 Although to deal with the discrete system is technically much more simple,
 the treatment of the discrete inverse problem is quite different to the continuous
 one and involve not so easy computations.
 Results have been obtained for the discrete Laplace equation which is a scalar one.
 Namely, for rectangular domains in two dimensions Oberlin proved uniqueness
 of the potential and gave a reconstruction algorithm in \cite{Ob}.
 Isozaki and Morioka generalized this result to $d$-dimension in \cite{Is-Mo}.

\section{Proof of the main result}

\subsection{Strategy of the proof}
 First, we write\refq{eq.HDu=lambdau} in the form\refq{eq20.H0u=-Vu}
 where the potential $V$ is a complex-valued invertible diagonal matrix function.
 Actually, it is easy to construct solutions of\refq{eq20.H0u=-Vu} by iterations
 on the first variable $n_1$.
 Then, we follow \cite{Is-Mo}.
 Specifically, by setting $p\in\Z$, we build a particular solution $v=(v^E,v^H)$
 to\refq{eq20.H0u=-Vu} which explicitly depends on the potential on the planes
 $n_1+n_2=p$, $n_1+n_2=p+1$, vanishes in the half-plane $n_1+n_2<p$,
 and does not depend on the potential on the boundary, except on $\partOm_1^+$.
 We prove that the input $v^{E/H}\tanOm$ only involves knowledge
 of the D-N operator and, finally, by iteration on $p$, that the output
 $\partanOm v^{E/H}$ provides knowledge of the potential.

\subsection{Treatment of the equations}
 We present the main equation\refq{eq.HDu=lambdau} in a more appropriate form.
 First, we write\refq{eq.HDu=lambdau} in the form $H_0u=Vu$ where the potential $V$
 is a complex-valued invertible diagonal matrix function.
 Then, using\refq{def.hatM} and\refq{def.H0u}, we write the linear system $H_0u=Vu$ as
 6 scalar equations with 6 unknown $u^E_j,u^H_j$, $1\le j\le 3$, where we break the symmetry
 between the three directions to make the direction $e_1$ preferred.
\begin{definition}
\label{def.admis}
 We say that a potential $V\in\calA(\Om;\calM_6(\CC))$ is admissible
 on $\Om$ when, for all $n\in \Om$, $V(n)$ is a diagonal $6\times 6$
 invertible matrix with complex coefficients. 
\end{definition}
 We write\refq{eq.HDu=lambdau} as
\begin{equation}
\label{eq1.H0u=-Vu}
\left\{ \begin{array}{rcl}
 M u^E &=& V^H u^H \\ M u^H &=& -V^E u^E
\end{array} \right. 
\end{equation}
 in $\Om$, where $V:=-\lambda(D_{\rm a})^{-1}|_\Om$ and $\lambda\neq 0$, so $V$ is admissible.
 We write $V=(V^E,V^H)$ with $V^{E/H}=$ diag $(V^{E/H}_1,V^{E/H}_2,V^{E/H}_3)$,
 and $\Lambda_V$ instead of $\Lambda(D_{\rm a})$. 
 
 In\refq{eq1.H0u=-Vu}, and in view of \refq{def.hatM}, the 1st and 4th scalar
 relations are treated specially, so we write\refq{eq1.H0u=-Vu} as:
\begin{equation}
\label{eq20.H0u=-Vu}
\left\{ \begin{split} (a) &
\left\{ \begin{array}{rcl}
 u^E_1   &=& (-2iV^E_1 )^{-1} (D_2 u^H_3 -  D_3 u^H_2)  \\
 u^H_1   &=& (2iV^H_1 )^{-1} (D_2 u^E_3 -  D_3 u^E_2) 
\end{array} \right.  \\
 (b)&
\left\{ \begin{array}{rcl}
 D_1 u^E_2  &=& 2iV^H_3 u^H_3 +  D_2 u^E_1  \\
 D_1 u^E_3 &=& - 2iV^H_2 u^H_2 +  D_3 u^E_1 \\
%
 D_1 u^H_2  &=& -2iV^E_3 u^E_3 +  D_2 u^H_1 \\
 D_1 u^H_3 &=& 2iV^E_2 u^E_2 +  D_3 u^H_1
\end{array} \right. 
\end{split}
\right.
\end{equation}
 in $\Om$.
\begin{remark}
\label{rem.uj}
 We observe that\refq{eq20.H0u=-Vu} in $\Om$ does not involve
 $u^{E/H}_j|_{\partOm_j}$, $j\in\{1,2,3\}$.
\end{remark}

\subsection{Modified Dirichlet-to-Neumann operator}
 Letting $\Gamma\subset\partOm$, we say that $f\in \calA(\Gamma;\CC^3)$
 is admissible on $\Gamma$ iff $f(n)\cdot\nu(n)=0$ for all $n\in \Gamma$.
 We denote by $\HadmE$ the set of such functions.
 Given $f=(f^E,f^H)\in (\HadmEp)^2$ and $g=(g^E,g^H)\in (\HadmEm)^2$,
 we consider the system with unknown $u$ that satisfies\refq{eq1.H0u=-Vu}
 in $\Om$ and the following partial boundary conditions:
\begin{equation}
\label{eq.BC2}
\left\{\begin{array}{lcll}
 u_\tanOm &=& f & {\rm on} \: \partOm\sauf\partOm_1^+\\
 \partanOm u &=& g & {\rm on} \: \partOm_1^-.
\end{array}\right.
\end{equation}
\begin{lemma}
\label{lem.22}
 System\refq{eq1.H0u=-Vu}-\refq{eq.BC2} admits a unique solution in $\calH(\Om;\CC^6)$.
\end{lemma}
 Proof. The above system is linear and square, so it suffices
 to prove that $f=0$ and $g=0$ imply $u=0$. We thus assume $(f,g)=0$.
 System\refq{eq20.H0u=-Vu} can be written
\begin{equation}
\label{eq21.a}
\left\{ \begin{array}{rcl}
 u^E_1(n) &=& (-2iV^E_1(n) )^{-1} \big(u^H_3(n+e_2)-u^H_3(n-e_2) -  u^H_2(n+e_3)+u^H_2(n-e_3)\big), \\
 u^H_1(n) &=& (2iV^H_1(n))^{-1} \big(u^E_3(n+e_2)-u^E_3(n-e_2) -  u^E_2(n+e_3)+u^E_2(n-e_3)\big) ,
\end{array} \right.
\end{equation}
 for $n\in \Om$ and
\begin{equation}
\label{eq21.b}
\left\{ \begin{array}{rcl}
 u^E_2(n+e_1) &=& 2iV^H_3 u^H_3(n) + u^E_1(n+e_2) - u^E_1(n-e_2) + u^E_2(n-e_1) ,  \\
 u^E_3(n+e_1) &=& -2iV^H_2 u^H_2(n) + u^E_1(n+e_3) - u^E_1(n-e_3) + u^E_3(n-e_1)  \\
 u^H_2(n+e_1) &=& -2iV^E_3 u^E_3(n) + u^H_1(n+e_2) - u^H_1(n-e_2) + u^H_2(n-e_1) , \\
 u^H_3(n+e_1) &=& 2iV^E_2 u^E_2(n) + u^H_1(n+e_3) - u^H_1(n-e_3) + u^H_3(n-e_1) ,
\end{array} \right.
\end{equation}
 for $n\in \Om-e_1$.
 Fixing $n=(n_1,n')$ with $n_1=1$ and $n'\in\om:=[[1,R_2]]\times [[1,R_3]]$,
 we obtain $u^{E/H}_1(n)=0$, from\refq{eq21.a}-\refq{eq.BC2}.
 Fixing $n_1=2$ and $n'\in\om$, we obtain $u^{E/H}_2(n)=u^{E/H}_3(n)=0$, 
 from\refq{eq21.b}-\refq{eq.BC2}. By iteration on $n_1\in [[1,R_1+1]]$,
 we thus obtain $u^{E/H}_j=0$ in $\Om'\sauf\partOm_j$, $j\in [[1,3]]$.
 Hence, $u=0$ in $\calH(\Om;\CC^6)$. 
\qed

 We thus define the operator:
$$
 \tilde\Lambda_V \: : \quad (\HadmEp)^2 \times (\HadmEm)^2 \ni
 (f,g) \mapsto u_\tanOm \in  (\Hadm)^2,
$$
 where $u$ is the unique solution of\refq{eq1.H0u=-Vu}-\refq{eq.BC2}.
 We call $f$ (respect., $g$) a partial Dirichlet (respect., Neumann)
 admissible data.
\begin{coro}
\label{coro.23}
 Given the mapping $\Lambda_V$, a partial Dirichlet admissible data $f$
 on $\partOm\sauf\partOm_1^+$ and a partial Neumann admissible data $g$
 on $\partOm_1^-$, then there exists a unique Dirichlet (admissible)
 data $\tilde f$ on $\partOm$ such that $\tilde f=f$ on
 $\partOm\sauf\partOm_1^+$ and $\Lambda_V \tilde f=g$ on $\partOm_1^-$.
\end{coro}
 Proof. The mapping $Q$: $\tilde f\mapsto (\tilde  f|_{\partOm\sauf\partOm_1^+},
 (\Lambda_V \tilde f)|_{\partOm_1^-})$
 is linear from $(\Hadm)^2$ into $(\HadmEp)^2\times(\HadmEm)^2$.
 The two spaces have the same finite dimension.
 Thanks to Lemma \ref{lem.22}, $Q$ is injective.
 The conclusion follows.
\qed

\medskip
 Remark. As in \cite[Corollary 8.4]{Is-Mo}, the above operator $Q^{-1}$
 can be computed from $\Lambda_V$.

\subsection{Construction of solutions}
 We extend the potential outside $\Om$ by setting $V(n)=I_6$ for $n\not \in \Om$.\\
 Let $c_1\in\Z$ and fix $u^{E/H}_2$ and $u^{E/H}_3$ on $\{c_1-1,c_1\}\times\Z^2$.
 As shown by the proof of Lemma \ref{lem.22}, Relations\refq{eq21.a} and\refq{eq21.b}
 provide a solution $u$ to\refq{eq1.H0u=-Vu} in $[[c_1,\infty)\times\Z^2$.
 Let us rewrite\refq{eq21.a}. Assuming $n_1\ge c_1+1$ and using\refq{eq21.b},
 we compute
\begin{eqnarray*}
 u^E_1(n) &=&  (-2iV^E_1(n) )^{-1} \sum_\pm \pm u^H_3(n\pm e_2) \mp u^H_2(n\pm e_3) \\
 &=& (-2iV^E_1(n) )^{-1}  \sum_\pm \pm 2iV^E_2(n\pm e_2-e_1) u^E_2(n\pm e_2-e_1) \\
 && \pm u^H_1(n+e_3\pm e_2-e_1) \mp u^H_1(n-e_3\pm e_2-e_1) \pm u^H_3(n\pm e_2-2e_1) \\
 && \pm 2iV^E_3(n\pm e_3-e_1) u^E_3(n\pm e_3-e_1) \mp u^H_1(n+e_2\pm e_3-e_1)\\
 &&  \pm u^H_1(n-e_2\pm e_3-e_1) \mp u^H_2(n\pm e_3-2e_1) .
\end{eqnarray*} 
 Hence, replacing $n$ by $n+e_1$ we have, for $n_1\ge c_1$,
\begin{equation}
\label{eq22.a}
\begin{array}{rcl}
 u^E_1(n+e_1)
 &=&(-2iV^E_1(n+e_1) )^{-1} \sum_\pm \pm 2iV^E_2(n\pm e_2) u^E_2(n\pm e_2)\\
 &&  \pm 2iV^E_3(n\pm e_3)  u^E_3(n\pm e_3) \pm u^H_3(n\pm e_2-e_1) \\
 && \mp u^H_2(n\pm e_3-e_1) .
\end{array}
\end{equation}
 Similarly, we obtain
\begin{equation}
\label{eq22.b}
\begin{array}{rcl}
 u^H_1(n+e_1)  &=&(2iV^H_1(n+e_1) )^{-1} \sum_\pm \mp 2iV^H_2(n\pm e_2) u^H_2(n\pm e_2)\\
 &&  \mp 2iV^H_3(n\pm e_3) u^H_3(n\pm e_3) \pm u^E_3(n\pm e_2-e_1) \\
 && \mp u^E_2(n\pm e_3-e_1) .
\end{array}
\end{equation}
 Conversely, we obtain\refq{eq21.a} (for $n_1\ge c_1+1$) from\refq{eq22.a}-\refq{eq22.b}-\refq{eq21.b}.
 Thus, by setting $u^{E/H}$ on $n_1=c_1$, we can replace System\refq{eq21.a}-\refq{eq21.b}
 by System\refq{eq21.b}-\refq{eq22.a}-\refq{eq22.b}.

\medskip

\subsection{Domains of influence and domains of determination}
 We put, for $n\in \Z^3$, the following cones:
$$ C^-(n):=\{m\in\Z^3;\:  m_1 \le n_1 - |n_2-m_2|+|n_3-m_3|\} , $$
$$ C^+(n):=\{m\in\Z^3;\:  m_1 \ge n_1 + |n_2-m_2|+|n_3-m_3|\} . $$
 Letting $p\in \Z$, we also introduce the following definitions.

\begin{enumerate}
\item 
$
  \calE^\pm(p):=\{n\in \Z^3;\; \pm(n_1+n_2) >\pm p\}, \quad \calE_0(p):= \{n\in \Z^3;\; n_1+n_2= p\}.
$
\item  The Dirichlet (respect., Neumann) boundary value of a function
 $\varphi$: $\Om'\to \CC^l$ at $n\in \partOm$ is $\varphi(n)$ (respect.,
 $\varphi(n)-\varphi(m_\Om(n))$, where $m_\Om(n)\in\Om$ and $n\sim m_\Om(n)$.
\item  Let $u\in \calH(\Om;\CC^6)$.
 The sequence of partial Dirichlet and Neumann boundary values $\calD(u;p)$
 consists in the following values:
\begin{itemize}
\item $u^{E/H}_i(m)$ for $m\in \calE^+(p-1)\cap(\partOm\sauf\partOm_i)$, $i=1,2,3$;
\item $u^{E/H}_i(m)$ for $m\in \calE^+(p-1)$ and $m\sim (\partOm\sauf\partOm_i)$, $i=1,2,3$.
\end{itemize}
\item The sequence of partial values $\calV_j(V;p)$ of the potential consists
 in the values of $V^{E/H}_j|_{\calE^+(p-1)\cap \Om}$, $j=1,2,3$.
\end{enumerate}
\begin{lemma}
\label{lem.23}
 Let $N\in \Z^3$ with $N_1\ge 1$ and $u$ satisfy\refq{eq21.b}-\refq{eq22.a}-\refq{eq22.b}
 in $[[0,N_1-1]]\times\Z^2$.
 Then, $u(N)$ depends only on $\{u^{E/H}_j(m)$; $j\neq 1$ $m\in C^-(N)\cap (\{-1,0\}\times\Z^2)\}$,
 $\{u^{E/H}_1(m)$; $m\in C^-(N) \cap(\{0\}\times\Z^2) \}$ and on $V|_{C^-(N)}$,
 and, in fact, $u(N)$ can be expressed as a linear combination of the $u^{E/H}_j(m)$,
 for $j\neq 1$ and $m\in C^-(N)\cap (\{-1,0\}\times\Z^2)$, and of the $u^{E/H}_1(N)$
 for $m\in C^-(N) \cap(\{0\}\times\Z^2)$, with coefficients depending only on $V|_{C^-(N)}$.\\
 In particular, if $u^{E/H}_2$ and $u^{E/H}_3$ vanish in $C^-(N)\cap (\{-1,0\}\times\Z^2)$,
 and if $u^{E/H}_1$ vanishes in $C^-(N)\cap (\{0\}\times\Z^2)$,
 then, $u=0$ in $C^+(N)\cap  ([[1,N_1]]\times\Z^2)$.
\end{lemma} 
 Proof.
 We see that $u(N)$ depends only on $\{u(m)$; $m\neq N$ and $m\in C^-(N)\}$.
 By a backward iteration on the first variable from $N_1$ to $1$,
 the result follows.

\begin{coro}
\label{coro.2.5}
 Let $p\in\Z$ and $u$ satisfy\refq{eq21.b}-\refq{eq22.a}-\refq{eq22.b} for $n_1\ge 0$.
 If $u^{E/H}_2$ vanishes in $\{-1\}\times (-\infty,p]]\times\Z$ and
 in $\{0\}\times (-\infty,p-1]]\times\Z$,
 if $u^{E/H}_3$ vanishes in $\{-1,0\}\times \Z^2$ and if $u^{E/H}_1$ vanishes
 in $\{0\}\times (-\infty,p-1]]\times\Z$, then, $u=0$ in $\calE^-(p)\cap [[0,\infty)\times \Z^2$.
\end{coro} 
 Similarly, we have
\begin{lemma}
\label{lem.2.6}
 Let $u$ be a solution of $H_0u=Vu$ in $\Om$. Let $n\in\Om$.\\
1) Assume $n_1=1$.

 a) The values $u^{E/H}_j(n)$ for $j\neq 1$ are known from the Dirichlet
 and Neumann boundary values $u^{E/H}_j$ at $n-e_1\in \partOm_1^-$.

 b) The value $u^{E/H}_1(n)$ depends only on $V^{E/H}_1(n)$,
 on the Dirichlet and Neumann boundary values of $u^{E/H}_j$ on
 $\partOm_1^-\cap C^-(n)$, $j\in \{2,3\}$, and on the Dirichlet boundary
 values of $u^{E/H}_k$ on $\partOm_j \cap C^-(n+e_1)\cap (\{1\}\times\Z^2)$,
 where $k\neq j$, $j,k\in \{2,3\}$.\\
2) Assume $1<n_1\le R_1$. 
 Then, $u(n)$ depends only on $V|_{\Om\cap C^-(n)}$, on the Dirichlet
 and Neumann boundary values of $u^{E/H}_j$ on $\partOm_1^-\cap C^-(n)$,
 $j\in \{2,3\}$, and on the Dirichlet boundary values of $u^{E/H}_k$
 at $\partOm_j\cap C^-(n)$, where $k\neq j$, $j\in \{2,3\}$.
\end{lemma}
 Proof.
 1) a) Obvious. Actually, we have, for example,
$$ u^{E/H}_2(1,n') = (u^{E/H}_\tanOm - \partanOm u^{E/H})_3(0,n') .$$

 b) It is directly the consequence of Formula\refq{eq21.a}.\\
 2) Except the case $n_1=2$, $n_2\in \{1,R_2\}$ or $n_3\in \{1,R_3\}$ for which
 $n\sim\partOm_2\cup\partOm_3$ so the value $u^{E/H}_1$ depends directly
 on the Dirichlet and Neumann boundary values $u^{E/H}_\tanOm(m)$ and $\partanOm u^{E/H}(m)$
 at $m\in\partOm_2\cup\partOm_3$ such that $n\sim m$, the result is easily verified
 by iterating over $n_1$ from $2$ to $R_1$ using Formulas\refq{eq21.b}-\refq{eq22.a}-\refq{eq22.b}.
\qed

\begin{lemma}
\label{lem.2.8}
 Let $u\in \calH(\Om,\CC^6)$ be a solution of $H_0u=Vu$ in $\Om$.\\
 Let $q\in\Z$.
 The function $u|_{\calE^+(q-1)\cap\Om}$ depends only on $\cup_j\calV_j(V;q-1)$,
 $\calD(u;q)$ and $u^{E/H}_3|_{\partOm_1^+\cap \calE_0(q-1)}$,
 and, in fact, for all $n\in \calE^+(q-1)\cap\Om$, $u(n)$ can be expressed
 as a linear combination of elements of $\calD(u;q)$ and of 
 $(u^{E/H}_3(m))_{m\in \partOm_1^+\cap \calE_0(q-1)}$ with coefficients
 depending only on $\cup_j\calV_j(V;q-1)$.
%
\end{lemma} 
 Proof.
 In \refq{eq21.b} we can exchange the vectors $e_1$ and $-e_1$, i.e., $u^{E/H}_j(n\pm e_1)$
 and $-u^{E/H}_j(n\mp e_1)$, $j\neq 1$. We obtain
\begin{equation}
\label{eq31}
\left\{ \begin{array}{rcl}
 u^E_2(n) &=& -2iV^H_3(n+e_1) u^H_3(n+e_1) - u^E_1(n+e_2+e_1) + u^E_1(n-e_2+e_1) - u^E_2(n+2e_1) ,  \\
 u^E_3(n) &=& 2iV^H_2(n+e_1) u^H_2(n+e_1) - u^E_1(n+e_3+e_1) + u^E_1(n-e_3+e_1) - u^E_3(n+2e_1),
\end{array} \right.
\end{equation}
 where $n_1\ge 0$; the same applies to $u^H_2$ and $u^H_3$.
 In addition, for $n_1<R_1$, we have relations similar to\refq{eq22.a} and\refq{eq22.b}
 showing that $u^E_1(n)$ and $u^H_1(n)$ depend only on $V^{E/H}_1(n)$,
 $V^{E/H}_j(n+e_1\pm e_j)$, $u^{E/H}_j(n+e_1\pm e_j)$ and on
 $u^{E/H}_j(n+2e_1\pm e_k)$ for $j\neq 1$, where $k\not\in\{1,j\}$;
 actually, we have, for $1\le n_1\le R_1-1$,
\begin{equation}
\label{eq32}
\begin{array}{rcl}
 u^E_1(n)
 &=& (-2iV^E_1(n) )^{-1} \sum_\pm \mp 2iV^E_2(n\pm e_2+e_1) u^E_2(n\pm e_2+e_1)\\
 &&  \mp 2iV^E_3(n\pm e_3+e_1)  u^E_3(n\pm e_3+e_1) \pm u^H_3(n\pm e_2+2e_1) \\
 && \mp u^H_2(n\pm e_3+2e_1).
\end{array}
\end{equation}
 (For $n_1=R_1$, we use\refq{eq21.a}.)
 Consider the version of Lemma \ref{lem.2.6} where the first direction $e_1$
 (respect., $C^-(N)$) is swapped with $-e_1$  (respect., $C^+(N)$).
 But we have $C^+(N)\subset \calE^+(q-1)$ for all $N\in \calE^+(q-1)$.
  The conclusion follows.
\qed

\subsection{Particular solutions}
 Let $\tau>0$ and $p\in [[1,R_1+R_2+1]]$.
 We consider the solution $u$ of\refq{eq21.b},\refq{eq22.a}\refq{eq22.b} for $n_1\ge0$
 and initialized by the following conditions:
\begin{equation}
\label{def.wp}
\left\{
\begin{array}{rcll}
 u_3^{E/H}(n) =  u_2^{E/H}(n) &=& 0 & n_1=-1,\\
 u_1^{H}(n) = -u_2^{H}(n)&=& i^{n_3}e^{\tau n_3}\delta_p(n_2) & n_1=0 ,\\
 u_1^{E}(n) = -u_2^{E}(n)&=& -i^{n_3}e^{\tau n_3}\delta_p(n_2) & n_1=0  ,\\
 u_3^{E/H}(n)&=& 0 & n_1=0.
\end{array} \right.
\end{equation}
 We set $w(;p)=u$ as a function of $\calH(\Om;\CC^6)$.
 It satisfies\refq{eq1.H0u=-Vu}.
%
%
%
\medskip
 Since $w^{E/H}(n;p)$ may depend on unknown values of $V|_{\Om}$ for
 $n\in \partOm_2^+$, and, in case, in addition, we want to impose homogeneous
 boundary conditions on $\partOm_3$, we consider the function $v(;p)\in \calH(\Om;\CC^6)$,
 which is, thanks to Lemma \ref{lem.22}, the unique solution of\refq{eq1.H0u=-Vu}
 in $\Om$ with the following boundary conditions:
\begin{equation}
\label{def.vp}
\left\{
\begin{array}{rcll}
 v^{E/H}_\tanOm(n;p) &=& w^{E/H}_\tanOm(n;p) & n\in \partOm_1^-,\\
 \partanOm v^{E/H}(n;p) &=& \partanOm w^{E/H}(n;p) & n\in \partOm_1^-,\\
 v^{E/H}_1(n;p) &=& w^{E/H}_1(n;p) & n\in \partOm_2^+\cap \calE_0(p),\\
 v^{E/H}_3(n;p) &=& w^{E/H}_3(n;p) & n\in \partOm_2^+\cap \calE_0(p+1),\\
 v^{E/H}_j(n;p) &=& 0 & \mbox{at other values $(n,j)\in \partOm\times[[1,3]]$}.
\end{array} \right.
\end{equation}
 Remark. We thus have
$$
 v^{E/H}_j(n;p) = w^{E/H}_j(n;p),\quad n \in\{0,1\}\times [[1,R_2]]\times[[1,R_3]],
  \quad{\rm and} \quad j\neq 1.
$$
 For simplicity, if there is no confusion, we write $v=v(;p)$.
\begin{lemma}
\label{lem.32}
 There exist $C^{E/H}(v)$ and $D^{E/H}(v)$ defined on $\calE_0(p+1)$ and $F^{E/H}(v)$ defined
 on $\calE_0(p)$ which depend only on $\calD(v;p+1)$ and $\calV(V,p+2)$,
 and, in fact, are linear combinations of the elements of $\calD(v;p+1)$ with coefficients
 depending only on $\calV(V,p+2)$, such that, for $n_1\ge 1$, we have
 the following relations:
\begin{equation}
\label{rel1.wp}
\left\{
\begin{array}{rcll}
  v(n) &=& 0 & n\in \calE^-(p),\\
  v^{E/H}_3(n) &=& 0 & n\in \calE_0(p);
\end{array} \right.
\end{equation}
\begin{equation}
\label{rel2.wp}
\left\{
\begin{array}{rcll}
 v^{E/H}_2(n) &=& v^{E/H}_1(n-e_1+e_2)  &  n\in \calE_0(p) ,\\
  v^{E/H}_1(n)   &=& -\frac{V^{E/H}_2(n-e_1+e_2)}{V^{E/H}_1(n)}v^{E/H}_2(n-e_1+e_2)
 &  n\in \calE_0(p);\\
\end{array} \right.
\end{equation}
\begin{equation}
\label{rel3.wp}
\left\{
\begin{array}{rcll}
 2iV_3^H(n)v^H_3(n) &=& v_1^E(n-e_2) - v_2^E(n-e_1) + C^E(v;n)  & n\in \calE_0(p+1) , \\
 -2iV_3^E(n)v^E_3(n) &=& v_1^H(n-e_2) - v_2^H(n-e_1) + C^H(v;n)  & n\in \calE_0(p+1) , \\
 v^{E/H}_3(n) &=& D^{E/H}(v;n)  & n\in \calE_0(p+1) ,
\end{array} \right.
\end{equation}
\begin{equation}
\label{rel4.wp}
\left\{
\begin{array}{rcll}
 -2iV_2^H(n-e_1)v_2^H(n-e_1) + v_1^E(n+e_3-e_1) - v_1^E(n-e_3-e_1)&=& v^E_3(n)  & n\in \calE_0(p+1),\\
  2iV_2^E(n-e_1)v_2^E(n-e_1) + v_1^H(n+e_3-e_1) - v_1^H(n-e_3-e_1)  &=& v^H_3(n) & n\in \calE_0(p+1);
\end{array} \right.
\end{equation}
\begin{equation}
\label{rel5.wp}
\left\{
\begin{array}{rcll}
 -2iV_2^{E}(n)v^{E}_2(n)  - v_2^{H}(n+e_3)  + v_2^{H}(n-e_3)  &=& F^{E}(v;n) & n\in \calE_0(p) ,\\
 2iV_2^{H}(n)v^{H}_2(n)  - v_2^{E}(n+e_3)  + v_2^{E}(n-e_3)  &=& F^{H}(v;n) & n\in \calE_0(p) .
\end{array} \right.
\end{equation}
\end{lemma}
 Proof. Put $u=v(;p)-w(;p)$.
 Thanks to Corollary \ref{coro.2.5} we have $w(n;p)=0$ for $n\in \calE^-(p)$.
 Let $n_1\ge0$ and $n\in \calE_0(p-1)$. The 2nd and 4th relations in\refq{eq21.b}
 and that $w(;p)|_{\calE^-(p)}=0$ imply $w(;p)^{E/H}_3(n+e_1)=0$.
 Hence, $w^{E/H}_3(;p) = 0$ in $\calE_0(p)$.
 Then, $u$ satisfies\refq{eq1.H0u=-Vu} in $\Om$ and
\begin{equation}
\label{eq.v-w}
\left\{
\begin{array}{rcll}
 u^{E/H}_\tanOm(n) &=& 0 & n\in \partOm_1^-,\\
 \partanOm u^{E/H}(n) &=& 0 & n\in \partOm_1^-,\\
 u^{E/H}_1(n) &=& 0 & n\in \partOm_2^+\cap \calE_0(p),\\
 u^{E/H}_1(n) &=& 0 & n\in \partOm_2\cap \calE^-(p-1),\\
 u^{E/H}_3(n) &=& 0 & n\in \partOm_2^+\cap \calE_0(p+1),\\
 u^{E/H}_3(n) &=& 0 & n\in \partOm_2\cap \calE^-(p).
\end{array} \right.
\end{equation}
 Thanks to Lemma \ref{lem.2.6} we obtain that $u$ vanishes in $\calE^-(p+1)$.
 Hence, it suffices to prove\refq{rel1.wp}$\sim$\refq{rel5.wp} with $v$
 replaced by $w(;p)$ and $\calD(v;p+1)$ replaced by $\calD(w(;p),p+1)$.
 We also extend $v_j^{E/H}$ by setting $v_j^{E/H}(n)=w_j^{E/H}(n;p)$ at the points
 $n\in \calE^+(p-1)\cap [[-1,0]]\times\Z^2$ where $w_j^{E/H}(;p)$ is already
 defined and not $v_j^{E/H}$.  
\begin{itemize}
\item 
 Consequently,\refq{rel1.wp} is proved since it is proved for $w(;p)$.
\item Let $n_1\ge1$ and $n\in \calE_0(p-1)$.
 The first and 3rd relations in\refq{eq21.b} and\refq{rel1.wp} imply
 the first relation in\refq{rel2.wp} with $n$ replaced by $n+e_2$.
 Let $n_1\ge0$ and $n\in \calE_0(p-1)$. Relations\refq{eq22.a},\refq{eq22.b}
 and\refq{rel1.wp} imply the second relation in\refq{rel2.wp} (with $n$
 replaced by $n+e_1$).
 Thus,\refq{rel2.wp} is proved.
\item Let $n_1\ge1$ and $n\in \calE_0(p+1)$. The 1st relation in\refq{eq21.b}
 and\refq{rel2.wp} implies:
\begin{eqnarray*}
 2iV_3^H(n)v^H_3(n) = v_1^E(n-e_2) - v_2^E(n-e_1) - v_1^E(n+e_2) + v_2^E(n+e_1) ,
\end{eqnarray*}
 which provides the 1st relations of\refq{rel3.wp} with
 $C^H(n):= - v_1^E(n+e_2) + v_2^E(n+e_1)$.
 The same applies to $2iV_3^E(n)v^E_3(n)$. Thanks to Lemma \ref{lem.2.8}, and since
 $n+e_1$ and $n+e_2$ belong to $\calE_0(p+2)$, we see that $C^{E/H}|_{\calE_0(p+1)}$
 depends only on $\calD(v;p+1)$ and $\calV(V,p+2)$.\\
 Let $n_1\ge2$ and $n\in \calE_0(p+2)$. The 2nd relation in\refq{eq21.b}  implies
\begin{eqnarray*}
 v_3^E(n-e_1) &=& v^E_3(n+e_1) + 2iV_2^H(n) v_2^H(n) - v_1^E(n+e_3) + v_1^E(n-e_3) \\
 &=:&  D^E(n-e_1).
\end{eqnarray*}
 The same applies to $v_3^H(n-e_1)$.
 We obtain the 3rd relation in\refq{rel3.wp} with $n$ replaced by $n-e_1\in \calE_0(p+1)$.
 Thus, $D^{E/H}|_{\calE_0(p+1)}$ depends only on $v_|{\calE^+(p+1)}$ and $V_2^H|_{\calE_0(p+2)}$.
 Thanks to Lemma \ref{lem.2.8}, we then see that $D^{E/H}|_{\calE_0(p+1)}$ depends
 only on $\calD(v;p+1)$ and $\calV(V,p+2)$.
 Thus,\refq{rel3.wp} is proved.
\item Let $n_1\ge0$ and $n\in \calE_0(p)$. The 2nd and 4th relations in\refq{eq21.b}
 and\refq{rel1.wp} imply the first relation in\refq{rel4.wp}.
 Similarly, with the exponents $E$ and $H$ swapped.
 Thus,\refq{rel4.wp} is proved (with $n$ replaced by $n+e_1$).
\item Let $n_1\ge0$ and $n\in \calE_0(p+1)$.
 Relations\refq{eq22.a}-\refq{rel1.wp}-\refq{rel2.wp} provide:
\begin{eqnarray*}
  -2iV_2^E(n-e_2)v^E_2(n-e_2)  - v_2^H(n+e_3-e_2) + v_2^H(n-e_3-e_2) \\
   + 2iV_3^E(n+e_3)v^E_3(n+e_3) - 2iV_3^E(n-e_3)v^E_3(n-e_3)  + v_3^H(n+e_2-e_1)  \\
  =  -2iV^E_1(n+e_1)v^E_1(n+e_1)  -2iV^E_2(n+e_2) v^E_2(n+e_2) 
  =: \tilde F^E(n-e_2).
\end{eqnarray*}
 We see that $\tilde F^E$ is defined on $\calE_0(p)$ and depends only on
 $v|_{\calE^+(p+1)}$ and $V|_{\calE_0(p+2)}$, and so, thanks to Lemma \ref{lem.2.8},
 only on $\calD(v;p+1)$ and $\calV(V,p+2)$.
 We replace in the above relation $2iV^{E/H}_3v^{E/H}_3$ and $v^{E/H}_3$
 by the corresponding expressions in\refq{rel3.wp}.
 We thus obtain
\begin{eqnarray*}
 -2iV_2^E(n-e_2)v^E_2(n-e_2)  - v_2^H(n+e_3-e_2)  + v_2^H(n-e_3-e_2) \\ \\
  = \tilde F^E(n-e_2)  - C^E(n+e_3) + C^E(n-e_3) - D^H(n+e_2-e_1) \\
  =: F^E(n-e_2).
\end{eqnarray*}
 Hence, $F^E$ is defined on $\calE_0(p)$ and depends only on $\calD(v;p+1)$
 and $\calV(V,p+2)$. Similarly, with the exponents $E$ and $H$ swapped.
 Replacing $n-e_2$ by $n\in \calE_0(p)$ in the above relation we obtain\refq{rel5.wp}.
 Thus,\refq{rel5.wp} is proved.
\end{itemize}
\qed
\begin{remark}
\label{rem.v=w}
 In Lemma \refq{lem.32} we can consider that $n$ is in $\Om$ or not
 since we can extend $v$ by setting $v(n)=w(n;p)$ for $n\in \calE_0(p)$
 and $v^{E/H}_3(n)=w^{E/H}_3(n;p)$ for $n\in \calE_0(p+1)$. 
\end{remark}

\medskip
 Letting $p\in \Z$, $(n_1,n_3)\in \Z^2$, $\tau\in\R$, we set:
\begin{equation}
\label{def.AEH}
  A^{E/H}(n_1,p,n_3) := \Pi_{j=-1}^{n_1} \frac{V^{E/H}_2(j-1,p-j+1,n_3)}{V^{E/H}_1(j,p-j,n_3)} .
\end{equation}
 (For $n_1\le -1$ we put $A^{E/H}(n_1,p,n_3)=1$.)
\begin{equation}
\label{def.Kt}
 K^{E/H}_\tau(n_1,p,n_3)= A^{E/H}(n_1,p,n_3+1) + e^{-2\tau} A^{E/H}(n_1,p,n_3-1).
\end{equation}
\begin{coro}
\label{coro.31}
 (1) For $n_1\ge 1$ we obtain the following relations.
\begin{equation}
\label{rel6.wp}
\left\{
\begin{array}{rcll}
 v^E_1(n) &=&  (-1)^{n_1-1} i^{n_3}e^{\tau n_3} A^E(n_1,p,n_3) & n\in \calE_0(p),\\
 v^E_2(n) &=&  (-1)^{n_1} i^{n_3}e^{\tau n_3} A^E(n_1-1,p,n_3) & n\in \calE_0(p),\\
 v^H_1(n) &=&  (-1)^{n_1} i^{n_3}e^{\tau n_3} A^H(n_1,p,n_3) & n\in \calE_0(p),\\
 v^H_2(n) &=& (-1)^{n_1-1} i^{n_3}e^{\tau n_3} A^H(n_1-1,p,n_3) & n\in \calE_0(p),
\end{array} \right.
\end{equation}
\begin{equation}
\label{val.V2}
 2V_2^E(n) A^E(n_1-1,p,n_3) = L^{E/H}(v;n;\tau) + e^{\tau} A^H(n_1-1,p,n_3+1)
   + e^{-\tau} A^H(n_1-1,p,n_3-1),
\end{equation}
%
%
%
\begin{equation}
\label{rel7.wp}
 K^{E/H}_\tau(n_1,p,n_3) + K^{E/H}_\tau(n_1-1,p,n_3) = T^{E/H}(v;n;\tau),
\end{equation}
 where $L^{E/H}(v)$ and $T^{E/H}(v)$, defined on $\calE_0(p)$, depend only
 on $\calD(v;p+1)$ and $\calV(V,p+2)$.\\
(2) For $\tau$ sufficiently large, we have $v^{E/H}_3(n)\neq 0$,
 for all $n\in \calE_0(p+1)\cap \Om$.
%
%
\end{coro}
 Proof. As above, it suffices to prove the relations with $v$ replaced by $w=w(;p)$.
 Let $n\in \calE_0(p)$.\\
(1)? a)
  Using\refq{rel2.wp} and\refq{def.wp}, we easily obtain\refq{rel6.wp}
  by iteration on $n_1$.\\
b)
 Using\refq{rel5.wp} and\refq{rel6.wp}, we have:
\begin{eqnarray*}
\nonumber  F^E(v;n) &=&   -2iV_2^E(n)v^E_2(n)  - v_2^H(n+e_3)  + v_2^H(n-e_3) \\
\label{rel.FE}
  &=& (-1)^{n_1} i^{n_3+1}e^{\tau n_3}\{-2V_2^E(n)( A^E(n_1-1,p,n_3) \\
\nonumber &&  + e^{\tau} A^H(n_1-1,p,n_3+1)  + e^{-\tau} A^H(n_1-1,p,n_3-1) \}  .
\end{eqnarray*}
 Putting 
\begin{eqnarray*}
 L^E(v;n;\tau) &:=& (-1)^{n_1+1} (-i)^{n_3+1}e^{-\tau n_3}F^E(v;n).
\end{eqnarray*}
 The same applies to $L^H$, so we get\refq{val.V2}.
 
c) Using\refq{rel5.wp},\refq{rel4.wp}, the third relation in\refq{rel3.wp}
 and\refq{rel2.wp}, we compute:
\begin{eqnarray*}
 F^E(v;n) &=&   -2iV_2^E(n)v^E_2(n)  - v_2^H(n+e_3)  + v_2^H(n-e_3) \\
 &=&  (v_1^H(n+e_3) - v_1^H(n-e_3)  - D^H(n+e_1)) - v_2^H(n+e_3)  + v_2^H(n-e_3)\\
 &=&  v_2^H(n+e_1-e_2+e_3) - v_1^H(n+e_1-e_2-e_3) - v_2^H(n+e_3) \\
 &&  + v_2^H(n-e_3) - D^H(v;n+e_1).
\end{eqnarray*}
 Hence, using\refq{rel6.wp},\refq{def.Kt}, we compute
\begin{eqnarray*}
  G^H(v;n) &:=& F^E(v;n) + D^H(v;n+e_1)\\
 &=&  v_2^H(n+e_1-e_2+e_3) - v_2^H(n+e_1-e_2-e_3) - v_2^H(n+e_3)  + v_2^H(n-e_3)  \\ 
 &=& (-1)^{n_1} i^{n_3+1}e^{\tau (n_3+1)} \{A^H(n_1,p,n_3+1)+A^H(n_1-1,p,n_3+1) \\
 && +  e^{-2\tau} \big(A^H(n_1,p,n_3-1)+A^H(n_1-1,p,n_3-1)\big)\}\\
 &=& (-1)^{n_1} i^{n_3+1}e^{\tau (n_3+1)} \big( K^H_\tau(n_1,p,n_3) + K^H_\tau(n_1-1,p,n_3) \big).
\end{eqnarray*}
 Putting
$$ T^H(v;n)=G^H(v;n)(-1)^{n_1} (-i)^{n_3+1}e^{-\tau(n_3+1)} .$$
 The same applies to $T^E$, $G^E$, and the conclusion follows.\\
(2)? Let $n_1\ge 1$ and $n_2=p+1-n_1$. Thanks to\refq{rel4.wp}
 and\refq{rel6.wp}, we have
\begin{eqnarray}
\nonumber
   v^E_3(n) &=& -2iV^H_2(n_1-1,n_2,n_3) v^H_2(n_1-1,n_2,n_3)\\
\nonumber
&&  + v^E_1(n_1-1,n_2,n_3+1)  - v^E_1(n_1-1,n_2,n_3-1) \\
&=& 
\label{val.uE3}
 (-1)^{n_1} i^{n_3+1}e^{\tau (n_3+1)} \Big\{A^E(n_1-1,p,n_3+1)
  \\
\nonumber
 && - e^{-\tau} 2V^H_2(n_1-1,n_2,n_3) A^H(n_1-2,p,n_3) \\
\nonumber
 && + e^{-2\tau} A^E(n_1-1,p,n_3-1)  \Big\}.
\end{eqnarray}
 The same applies to $v^H_3(n)$.
 We thus see that, for $\tau>\tau_0$ where $\tau_0$ depends only
 on $\|V\|_\infty$ and on the size of $\Om$, $v^{E/H}_3(n)\neq 0$. 
\qed

\begin{lemma}
\label{lem.31}
 The Dirichlet and Neumann values of $v=v(;p)$ on $\partOm\sauf\partOm_1^-$
 and the Dirichlet values of $v$ on $\partOm_2\cup\partOm_3$
 do not depend on $V|_{\Om}$.
\end{lemma}
 Proof.
1)
 We have $v^E_3=0$ on $\partOm_1^-$ and, for $n_1=1$, Relation\refq{eq21.b} provides
\begin{eqnarray*}
  v^E_3(n) &=&  -2iV_2^H v^H_2(n-e_1) = 2(-i)^{n_3+1}e^{\tau n_3} \delta_p(n_2).
\end{eqnarray*}
 The same applies to $v^H_3$.
 Hence, the Dirichlet and Neumann boundary values of $v^{E/H}_3$ on $\partOm_1^-$
 do not depend on $V$.\\

 The values $v^{E/H}_2|_{n_1=0}$ are set independently of $V|_{\Om}$,
 and Relation\refq{rel2.wp} provides for $n_1=1$:
\begin{eqnarray*}
  v^E_2(n) &=&  v^E_1(0,n_2+1,n_3) = -i^{n_3}e^{\tau n_3} \delta_{p-1}(n_2).
\end{eqnarray*}
 Similarly, $v^H_2(n) = i^{n_3}e^{\tau n_3} \delta_{p-1}(n_2)$.
 Hence, the values $v^{E/H}_j(n)$ and $(\partanOm v^{E/H})_j(n)$, for $n\in \partOm_1^-$
 and $j\neq 1$ do not depend on $V$.\\
2) Let $n\in \partOm_2^+$. If $n\not\in \calE^+(p-1)$ then $v^{E/H}_1(n)=0$.
 If $n\in \calE_0(p)$ then, thanks to \refq{rel6.wp}, we have
 $v^E_1(n)= (-1)^{n_1-1} i^{n_3}e^{\tau n_3}$ which is known,
 and if $n\in \calE^+(p)$ then $v^{E/H}_1(n)=0$ by construction.
 The same applies to $v^H_1$.\\
 If $n\not\in \calE^+(p)$ then $v^{E/H}_3(n)=0$.
 If $n\in \calE_0(p+1)$ then, thanks to \refq{val.uE3} and to the fact that
 $A^{E/H}(m_1,p,m_3)=1$ for $m_1\le n_1-1$, the value of
 $v^{E/H}_3(n)$ does not depend on $V$.
 If $n\in \calE^+(p+1)$ then $v^{E/H}_3(n)=0$ by construction.
 Hence, the values $v^{E/H}_j(n)$ for $n\in \partOm_2^+$ and $j\neq 2$ do
 not depend on $V$.\\
3) By construction, we have $v^{E/H}_j(n)=0$ for $n\in \partOm_2^-$,
 $j\neq 2$, and, $v^{E/H}_j(n)=0$ for $n\in \partOm_3$, $j\neq 3$.

 The lemma is proved.
\qed

\subsection*{Reconstruction of the potential}
 We make the reconstruction of $V$ by a backward iteration on $p$ from
 $R_1+R_2$ until $1$.
 We first observe that $\calV_j(V;q)$, $j\in [[1,3]]$, $q\ge R_1+R_2+1$, are known.
 Fix $p\in [[1,R_1+R_2]]$ and assume that $\calV_j(V;p+1)$, $j\neq3$,
 and $\calV_3(V;p+2)$ are known.
 Let us show how to compute $\calV_j(V;p)$, $j\neq 3$, and $\calV_3(V;p+1)$.
 It suffices to compute $V^{E/H}_j$, $j\neq 3$, on $\calE_0(p)$, and $V^{E/H}_3$
 on $\calE_0(p+1)$.
 We do it in two steps. We put $u:=v(;p)$. \\
1)  Take $n\in \calE_0(p)$ such that $n_1\ge 1$. \\
 a) By assumption on $\calV_j(V;p+1)$ and to the property on $T^{E/H}$
 in Corollary \ref{coro.31}, the left-hand side term of\refq{rel7.wp} can then be computed.
 Since $K^{E/H}_\tau(0,p,n_3)=1+ e^{-2\tau}$ we thus can compute
 $K^{E/H}_\tau(n_1,p,n_3)$ for $n_1\in [[1,R_1]]$ by iteration on $n_1$.
 Taking two different values of $\tau$ we see that we can compute, from\refq{def.Kt},
 $A^{E/H}(n_1,p,n_3)$ for $n_1\ge 1$ and $n_3\in [[1,R_3]]$.
 By assumption on $\calV_j(V;p+1)$ and to the property on $L^{E/H}$ in Corollary
 \ref{coro.31}, the right-hand side term of\refq{val.V2} can then be computed.
 Hence, $V^{E/H}_2$ can be computed on $\calE_0(p)$.\\
b) Thanks to the definition of $A^{E/H}$ in\refq{def.AEH}, we have
$$
 V^{E/H}_1(n_1,p-n_1,n_3) =
    \frac1{A^{E/H}(n_1,p,n_3)}
  \Pi_{j=0}^{n_1-1} \frac{V^{E/H}_2(j,p-j,n_3)}{V^{E/H}_1(j,p-j,n_3)}  .
$$
 Since, in addition, $V^{E/H}_1(0,p,n_3)=1$, $V^{E/H}_1$ can be computed on $\calE_0(p)$
 by iteration on $n_1$.

2) Let $n\in \calE_0(p+1)$ such that $n_1\ge 1$.
 The second and third relations of\refq{rel3.wp} and the fact that $u_3^E(n)\neq 0$
 provide:
\begin{equation}
 -2iV_3^E(n) = \frac1{D^E(u;n)} \big(u_1^H(n-e_2) - u_2^H(n-e_1) + C^H(u;n)\big) .
\end{equation}
 In addition, thanks to Corollary \ref{coro.31} and to the above step 1, the values
 $u_1^H(n-e_2)$ and $u_2^H(n-e_1)$ are computed from $V^H_j|_{\calE_0(p)}$, $j\neq 3$.
 Hence, $V_3^E(n)$ can be computed. Similarly,  $V_3^H(n)$ can be computed.
 The proof is achieved. \qed

\bibliographystyle{iopart_BibTeX_copie/iopart-num}

\end{document}